\documentclass{amsart}

\usepackage[colorlinks=true, urlcolor=blue,bookmarks=true,bookmarksopen=true,citecolor=blue,hypertex]{hyperref}
\usepackage{graphicx}
\usepackage{enumerate}
\usepackage{parskip}

\usepackage{amscd}
\usepackage{amssymb}
\usepackage{amsxtra}
\usepackage{amsmath}
\usepackage{enumerate}
\usepackage{mathrsfs}
\usepackage{centernot}

\usepackage{color}

\usepackage{amsmath,amsfonts,amssymb}
\usepackage{verbatim}
\usepackage{epsfig}
\usepackage{amsthm}
\usepackage{bbold}

\usepackage[T1]{fontenc}

\theoremstyle{plain}

\newtheorem{thm}{Theorem}

\newtheorem{cor}[thm]{Corollary}

\newtheorem{prop}[thm]{Proposition}

\theoremstyle{definition}

\theoremstyle{remark}

\newtheorem{rmk}[thm]{Remark}

\newtheorem{qs}[thm]{Question}

\theoremstyle{plain}

\newcommand{\Qed}{\hfill \qedsymbol \medskip}

\newcommand{\Id}{{{\mathchoice {\rm 1\mskip-4mu l} {\rm 1\mskip-4mu l}
      {\rm 1\mskip-4.5mu l} {\rm 1\mskip-5mu l}}}}
\newcommand\supp{\operatorname{supp}}

\newcommand\dist{\operatorname{dist}}

\def\ZZ{\mathbb{Z}}

\def\RR{\mathbb{R}}

\def\Id{\mathbb{1}}

\def\dd{\mathrm{d}}
\DeclareMathOperator{\Ham}{Ham}

\DeclareMathOperator{\Fix}{Fix}

\def\TT{{T^2}}
\def\TTs{\mathbb{T}}
\def\Ac{\mathcal{A}}

\DeclareMathOperator{\Area}{Area}

\begin{document}

\pagestyle{headings}

\bibliographystyle{alphanum}

\title[Hofer's distance between eggbeaters and autonomous Hamiltonians]{Hofer's distance between eggbeaters and autonomous Hamiltonian diffeomorphisms on surfaces}

\date{\today}

\author{Michael Khanevsky}
\address{Michael Khanevsky, Faculty of Mathematics,
Technion - Israel Institute of Technology
Haifa, Israel}
\email{khanev@technion.ac.il}

\begin{abstract}
Let $\Sigma$ be a compact surface of genus $g \geq 1$ equipped with an area form. We construct eggbeater Hamiltonian diffeomorphisms
which lie arbitrarily far in the Hofer metric from the set of autonomous Hamiltonians. This result is already known for $g \geq 2$
(our argument provides an alternative, very simple construction compared to previous publications) while the case $g = 1$ is new.
\end{abstract}

\maketitle

\section{Introduction}

Let $\Sigma$ be a compact surface of genus $g$, possibly with boundary, equipped with an area form $\omega$.
Denote by $\Ham(\Sigma)$ the group of Hamiltonian diffeomorphisms of $(\Sigma, \omega)$ that are supported in the interior of $\Sigma$.
For every smooth function $H\colon\Sigma\to\RR$ supported in the interior, there exists a unique vector field $X_H$ which satisfies
$$
dH(\cdot)=\omega(X_H,\cdot).
$$
It is easy to see that $X_H$ is tangent to the level sets of $H$. Let $h$ be the
time-one map of the flow $h_t$ generated by $X_H$. The diffeomorphism $h$ is
area-preserving, it belongs to $\Ham(\Sigma)$ and every Hamiltonian diffeomorphism arising this way is called
{\em autonomous}. Such a diffeomorphism is easy to understand
in terms of its generating Hamiltonian function.
We will denote the set of autonomous Hamiltonian diffeomorphisms by $\Ham_{Aut} (\Sigma)$.

Hofer's metric is a bi-invariant metric on $\Ham(\Sigma)$ (see the definition in Section~\ref{S:Hofer-Floer})
which roughly measures the `mechanical energy' needed to deform one Hamiltoian diffeomorphism into another.
Following \cite{Po-Sh:egg-beater}, we define
\[
  \text{Aut} (\Sigma) = \sup_{g \in \Ham(\Sigma)} \dist_H (g, \Ham_{Aut} (\Sigma))
\]
and address the following question:
\begin{qs} \label{Q:aut_dist}
  Is $\text{Aut} (\Sigma) < \infty$?
\end{qs}
Conjecture~1.1 of \cite{Po-Sh:egg-beater} asserts that the answer is negative.
Indeed, Question~\ref{Q:aut_dist} was negatively answered for most surfaces: firstly, in \cite{Po-Sh:egg-beater} for surfaces of genus $g \geq 4$,
later \cite{Ch:eggbeater} extended the argument to $g \geq 2$. In \cite{Ch-Me:Hof-entr} a similar construction to \cite{Po-Sh:egg-beater}
is analyzed from a different viewpoint, which also implies a proof for $g \geq 2$.
The case of an annulus is discussed in \cite{Bra-Kha:C0-gap}.

This text proposes a very simple argument which applies to $g \geq 1$:
\begin{thm} \label{T:main}
 Let $\Sigma$ be a compact symplectic surface of genus $g \geq 1$. Then for all $R > 0$ there exists $g_R \in \Ham (\Sigma)$
 such that the Hofer distance $\dist_H (g_R, \Ham_{Aut} (\Sigma)) > R$.
\end{thm}
Methods used in the publications mentioned above hold if one adds punctures to $\Sigma$,
hence the only two compact surfaces where the question remains open are the disk and the sphere.

We would like to add that previous publications answer Question~\ref{Q:aut_dist} as a part of much stronger statements.
Our argument addresses the question only in its basic version, which is in some measure compensated by simplicity of the proof.
We observe that if $g \in \Ham(\Sigma)$ has two fixed points whose trajectories must intersect for topological reasons,
such $g$ is non-autonomous. Then we construct an eggbeater map $g_R$ which features such fixed points. We show that
these fixed points of $g_R$ persist under
perturbations bounded by $R$ in the Hofer metric, hence the Hofer ball of radius $R$ around $g_R$ contains no autonomous maps.

\section{Recollection from the Hofer metric and Floer theory}\label{S:Hofer-Floer}

Let $(M, \omega)$ be a symplectic manifold, $g$ a Hamiltonian diffeomorphism with compact support in $M$.
The \emph{Hofer norm} $\|g\|$ (see ~\cite{Hof:TopProp}) is defined by
\[
  \|g\| = \inf_G \int_0^1 \max \left(G(\cdot, t) - \min G(\cdot, t) \right) \dd t
\]
where the infimum goes over all compactly supported Hamiltonian functions $G : M \times [0,1] \to \RR$
such that $g$ is the time-1 map of the induced flow. The Hofer distance is given by
\[
  \dist_H (g_1, g_2) = \| g_1 g_2^{-1}\|.
\]
and, intuitively, it expresses the amount of mechanical energy needed to deform the dynamics of $g_1$ into $g_2$.

\medskip

We remind the most basic ingredients of the Floer theory with a focus at non-contractible orbits.
Let $(\Sigma, \omega)$ be a compact symplectic surface of genus $g > 0$, for simplicity of notations we assume that
$\Area (\Sigma) = \int_\Sigma \omega = 1$. Pick $H_t : \Sigma \to \RR$
a time-dependent Hamiltonian function supported in the interior of $\Sigma$. We assume that $H_t$ is normalized:
$H_t$ has zero mean for all $t$ if $\Sigma$ is closed.
(For non-closed surfaces the normalization condition requires $H_t$ to vanish near the boundary and it automatically follows from our
assumptions on the support of $H_t$.)

Denote $S^1 = \RR/\ZZ$. We fix a non-contractible free homotopy class $\alpha \in \pi_1 (\Sigma)$ and its representative
$\gamma_\alpha : S^1 \to \Sigma$.
Let $\gamma : S^1 \to \Sigma$ be a loop with $[\gamma] = \alpha \in \pi_1 (M)$. Pick a smooth cylinder $w_\gamma : [0, 1] \times S^1  \to \Sigma$
with $w_\gamma(0, t) = \gamma (t)$ and $w_\gamma(1, t) = \gamma_\alpha (t)$ ($w_\gamma$ is called a \emph{cylindrical capping} of $\gamma$).
The Hamiltonian action of $(\gamma, w_\gamma)$ is defined by
\[
   \Ac_H(\gamma, w_\gamma) = \int_{[0, 1] \times S^1} (w_\gamma)^* \omega + \int_0^1 H_t (\gamma(t)) \dd t .
\]
By Stokes theorem, the action depends only on the homotopy class of $w_\gamma$ (as long as its ends $\gamma_\alpha, \gamma$ are fixed).
Hence for surfaces $\Sigma$ different from the torus $T^2 = S^1 \times S^1$ (namely, if $g > 1$ or if $g = 1$ and $\partial \Sigma \neq \emptyset$)
the action $\Ac(\gamma, w_\gamma)$ is independent of the capping $w_\gamma$ as there is a unique homotopy class of cylinders with fixed ends.
In the case $\Sigma = \TT$ there is a countable set of different homotopy classes and the action depends on a capping.
Changing a capping results in a shift of the action by an integer (a multiple of $\Area (\TT) = 1$).
We note that homotopy classes of cappings can be distinguished by looking at the local degree $\deg(w_\gamma, p)$ computed at a point
$p \in \TT \setminus (\gamma \cup \gamma_\alpha$).

Pick a generic $C^\infty$-smooth $S^1$-family $J_t$ of compatible almost complex structures on $\Sigma$. If $\partial \Sigma \neq \emptyset$
one asks that $J_t$ coincides near $\partial \Sigma$ with the standard cylindrical complex structure.
Let $u : \RR \times S^1 \to \Sigma$ be a cylinder. The Floer equation for $(H_t, J_t)$ and $u$ is
\[
  \mathcal{F}_{H,J}(u) = \partial_s u (s,t) + J_t (u(s,t)) (\partial_t u (s,t) - X_H (t, u(s,t)))  = 0
\]
and its solutions are called Floer cylinders. We note that in an open subset where $H_t (p)$ is constant in $p$, $X_{H_t} (p) = 0$ hence the Floer equation becomes the equation of a $J$-holomorphic curve.

Suppose that $H$ is non-degenerate. The energy of a Floer cylinder $u$ is given by
\[
  E(u) = \int_{\RR \times S^1} |\partial_s u|^2 \dd t \dd s =  \int_{\RR \times S^1} \omega \left(\partial_s u, J_t (\partial_s u)\right) \dd t \dd s .
\]
$E (u) < \infty$ implies that $u$ converges on either end to a $1$-periodic trajectory of the flow of $H_t$ (see \cite{Sa:HF-lectures}): there exist trajectories
$\gamma_+, \gamma_-$ of $X_{H_t}$ such that $\lim_{s \to \pm \infty} u (s, t) = \gamma_\pm (t)$. Moreover, the energy is equal to the action difference at the ends:
$$E(u) = \Ac (\gamma_+, w_{\gamma_-} \# u) - \Ac (\gamma_-, w_{\gamma_-})$$ where $w_{\gamma_-}$ is a capping of $\gamma_-$ and $w_{\gamma_-} \# u$
is a capping for $\gamma_+$ obtained by gluing $w_{\gamma_-}$ and $u$.
$E(u) \geq 0$ by its definition. $E(u) = 0$ if and only if the Floer cylinder $u (s, t)$ is trivial
(that is, $u (s,t) = \gamma_-(t) = \gamma_+(t)$ for all $s$).

Independence of the Hamiltonian isotopy: let $g \in \Ham (\Sigma)$. For surfaces $\Sigma \centernot\simeq S^2$, the group $\Ham (\Sigma)$ is simply connected
(see Chapter~7.2 in \cite{Po:geo-symp}). Therefore by Arnold conjecture for each $p \in \Fix (g)$, its trajectory $\gamma_{g_t,p} (t) = g_t (p)$
defines a free homotopy class $[\gamma_{g_t,p}] \in \pi_1 (\Sigma)$ which does not depend on the isotopy $g_t$ that connects $g$ to $\Id_\Sigma$.
By \cite{Sc:action-spectrum}, assuming that the Hamiltonian function generating the isotopy $g_t$ is normalized,
the action of $\gamma_{g_t,p}$ does not depend on $g_t$ as well
but only on $p \in \Fix (g)$ and, in the case of $\Sigma = \TT$, on the choice of a homotopy class of the capping.

\medskip

\cite{Ush:b-depth-appl} introduces the notion of \emph{boundary depth} for a Floer complex $CF(H, J)$ as
\[
 b (CF(H, J)) = \inf_{\beta \geq 0} \left\{ \forall \lambda \in \RR .
 CF^\lambda (H, J) \cap \partial_{H, J} (CF (H, J)) \subseteq \partial_{H, J} (CF^{\lambda+\beta} (H, J)) \right\} .
\]
Namely, it is the minimal action margin which is required to present every boundary chain in $CF (H, J)$ as a value of $\partial_{H, J}$.
$H_t$ is assumed to be non-degenerate, $J_t$ regular and, while \cite{Ush:b-depth-appl} considers contractible trajectories, the same
results apply to the Floer complex $CF_\alpha (H, J)$ for a non-contractible class of trajectories $\alpha \in \pi_1 (\Sigma)$.
In the case $\partial_{\alpha, H, J} \neq 0$, it is a straightforward corollary from the definition that
\[
 b (CF_\alpha(H, J)) \geq \inf \{E (u)\}
\]
where $u$ runs over all non-trivial Floer cylinders involved in the construction of $\partial_{\alpha, H, J}$.

The boundary depth $b (CF(H, J))$ is independent of the choice of $J_t$.
In the case of simply connected $\Ham(\Sigma)$, the value $b (CF(H, J))$ depends only on the time-$1$ map and not on the Hamitonian isotopy.
Moreover, by \cite{Ush:b-depth-appl} the boundary depth is $1$-Lipshitz with respect to the Hofer metric:
\[
	\left| b (CF_\alpha(H^1, J^1)) - b (CF_\alpha(H^2, J^2)) \right| \leq \dist_H (h_1, h_2)
\]
where $h_1, h_2 \in \Ham(\Sigma)$ are the respective time-$1$ maps for $H^1_t, H^2_t$.

\begin{cor} \label{Cor:stability}
  Let $\alpha \in \pi_1 (\Sigma)$ be a non-contractible homotopy class,
  $G_t : \Sigma \to \RR$ a Hamiltonian function which generates $g \in \Ham (\Sigma)$ and $J_t$ a generic compatible almost complex structure.
  Suppose that the set $\Fix_\alpha(g)$ of fixed points in the class $\alpha$ of $g$ is non-empty and these points are non-degenerate.
  Let $\epsilon = \inf \{E (u)\}$ where $u$ runs over all non-trivial Floer cylinders for $(G_t, J_t)$ which connect orbits of $\Fix_\alpha(g)$.
  Then any $h \in \Ham (\Sigma)$ with $\dist_H (g, h) < \epsilon$ has a fixed point in the class $\alpha$.
\end{cor}
Proof:
  the Floer homology $H(CF_\alpha (G, J)) = 0$ since it is independent of the pair $(G, J)$ and can be computed from $CF_\alpha (\Id_\Sigma, J')$
  which has no orbits in the class $\alpha$.
  Hence, assuming $\Fix_\alpha(g) \neq \emptyset$, the Floer differential $\partial_{\alpha, G, J} \neq 0$ and $b (CF_\alpha(G, J)) \geq \epsilon$.
  Pick $h \in \Ham (\Sigma)$ in the $\epsilon$-Hofer-neighborhood of $g$.
  If $h$ has a degenerate fixed point in $\Fix_\alpha(h)$, there is nothing to prove.
  Otherwise, selecting appropriate Floer data, $h$ induces a Floer complex for the class $\alpha$ with positive boundary depth.
  In particular, $h$ has a fixed point in the class $\alpha$.
\Qed

\section{Eggbeaters and autonomous maps}\label{S:eggs}

\subsection{Obstruction for autonomous diffeomorphisms}

Our argument is based on the following observation: Let $\alpha, \beta \in \pi_1 (\Sigma)$ be different free homotopy classes with
topologically non-trivial intersection. We mean by that that $\alpha, \beta$ do not admit disjoint representatives (for example,
this is the case when the intersection product $[\alpha] \cdot [\beta] \neq 0$). Then:

\begin{prop} \label{P:obstruction}
 Suppose $\alpha, \beta \in \pi_1 (\Sigma)$ have topologicallly non-trivial intersection.
 Let $g \in \Ham (\Sigma)$ be a diffeomorphism with fixed points $p_\alpha, p_\beta \in \Fix (g)$ whose trajectories represent
 $\alpha, \beta$, respectively. Then $g$ is not autonomous.
\end{prop}

Proof:
  suppose by contradiction that $g$ is generated by an autonomous flow $g_t$. Note that fixed points of autonomous flows arise in two possible scenarios:
  \emph{equilibrium} points which are stationary under the flow, or a point $p$ whose orbit under the flow is a simple loop $c$ in $\Sigma$.
  In the latter case the trajectory $\gamma_p (t) = g_t (p), \; t \in [0, 1] $ traverses $c$ once or several times.
  Moreover, all other points in $c$ have the same orbit (namely, $c$) and the same period as $p$.
  Hence $c \subseteq \Fix(g)$ and the trajectory of any $q \in c$ represents the same class $[\gamma_p] \in \pi_1 (\Sigma)$.

  $p_\alpha, p_\beta \in \Fix (g)$ are periodic points whose trajectories $\gamma_{p_\alpha}, \gamma_{p_\beta}$ represent
  $\alpha, \beta$, respectively. $\alpha, \beta$ have non-trivial intersection, hence both classes are non-contractible.
  In particular, $p_\alpha, p_\beta$ cannot be equilibrium points. Their orbits $c_{p_\alpha}, c_{p_\beta}$ are loops which contain periodic points
  from different homotopy classes, hence $c_{p_\alpha}, c_{p_\beta}$ are disjoint.
  Therefore $\gamma_{p_\alpha}, \gamma_{p_\beta}$ are disjoint as well, a contradiction.
\Qed

\subsection{An eggbeater map}
We start with the torus $\TTs = \RR^2 / \ZZ^2$ equipped with the symplectic form $\omega = \dd x \wedge \dd y$ so that $\Area (\TTs) = 1$.

Let $h : \RR / \ZZ \to [-1, 1]$ be a $C^\infty$-smooth function which satisfies the following:
\[
  h (t) = \left\{\begin{aligned}
	 1 - 50t^2, &\qquad t \in \left[ -0.01, 0.01\right] \\
	 1.25 - 5t, &\qquad t \in \left[ 1/8, 3/8\right] \\
	 -1 + 50(t-0.5)^2, &\qquad t \in \left[ 0.49, 0.51\right] \\
	 -3.75 + 5t, &\qquad t \in \left[ 5/8, 7/8\right]
  \end{aligned} \right.
\]
$h (t)$ is smoothly interpolated in the intervals
$t \in \left[ 0.01, 1/8\right] \cup \left[ 3/8, 0.49 \right] \cup \left[ 0.51, 5/8\right] \cup \left[ 7/8, 0.99 \right]$
so that the slope $h'(t)$ is strongly increasing for $t \in \left[ -1/8, 1/8\right]$ and strongly decreasing for $t \in \left[ 3/8, 5/8\right]$.
Example graph of such a function and its derivative $h'(t)$ is provided in the figure below.
We define Hamiltonian functions $F, P: \TTs \to \RR$ by $F(x, y) = h(y)$ and $P (x, y) =  h(x)$.
The autonomous flow induced by $F$ rotates horizontal fibers $\{y = c\}$ right/left with
anglular velocity $h'(c) \in [-5, 5]$ while the flow of $P$ rotates vertical fibers with velocity $-h'(c)$).
\begin{figure}[!htbp]
\begin{center}
\includegraphics[width=0.7\textwidth]{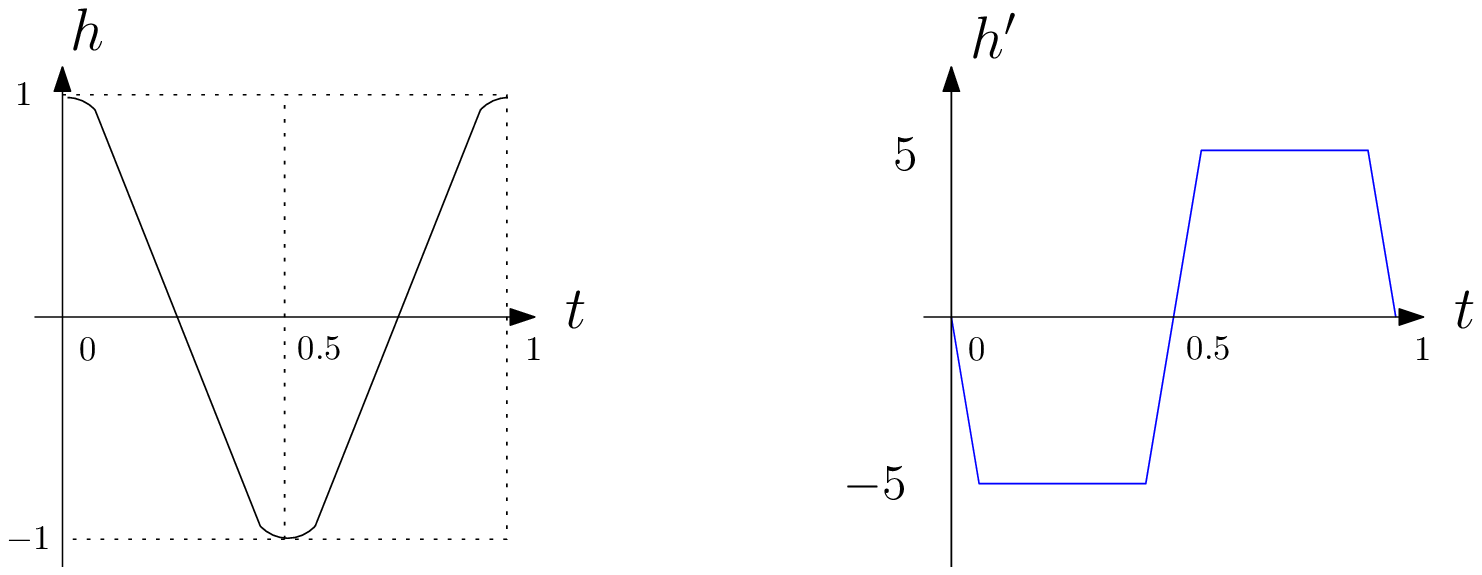}
\end{center}
\end{figure}

Pick a parameter $A > 1$, later on it will determine a lower bound for the distance to $\Ham_{Aut}$.
We construct an eggbeater map $g_A \in \Ham (\TTs)$ as follows.
Pick $r_A < \frac{1}{1000A}$ and let $q_0 = (0, 0)$. Note that $q_0$ is an equilibrium point for the flows induced by $F$ and $P$,
and both flows have velocity bounded by $\frac{1}{5A}$ in the disk $D' = \{ \|p - q_0\| < 2 r_A\}$.
We perturb the functions $F, P$ in $D'$ in a way that
they both vanish in the disk $D_A  = \{ \|p - q_0\| < r_A\}$ and remain normalized. That is, $\int_\TTs F \omega = \int_\TTs P \omega = 0$.
Denote by $\Phi_t$ and $\Psi_t$ the flows induced by the perturbed functions $F, P$ and
let $g_A = \Psi_{2A} \circ \Phi_A$ be a variant of the eggbeater map.
The perturbation is just a technical trick which is intended to make the argument below more elementary.
We note that the induced flows $\Phi_t, \Psi_t$ may have somewhat complicated dynamics in $D'$.
However, since the angular velocity $h'(c)$ of the unperturbed flows was very small in this neighborhood,
the perturbation did neither affect existing non-contractible $1$-periodic trajectories of $\Psi_{2A} \circ \Phi_A$ nor
created new ones.

\subsection{Surfaces of positive genus different from the torus}
Let $\Sigma$ be a compact symplectic surface of genus $g > 0$. We fix a parameter $A > 1$ and
adapt the construction above of the eggbeater map $g_A$ to $\Sigma$.

Recall that after the perturbation both functions $F, P$ and their flows are supported outside the disk $D_A$ around $q_0$.
In the case $g > 1$ we model $\Sigma$ by gluing a handlebody to $\TTs \setminus D_A$. Otherwise
(namely, $g = 1$ and $\Sigma$ is a torus with `holes') we cut small disks from $D_A$.
By picking an appropriate symplectic form in the complement of $\TTs \setminus D_A$
one ensures that the total area is $1$. The resulting surface $\widetilde{\Sigma}$ is symplectomorphic to $\Sigma$.
We restrict the perturbed functions $F, P : \TTs \to \RR$ constructed above to $\TTs \setminus D_A$ and then extend to
$\widetilde{F}, \widetilde{P} : \widetilde{\Sigma} \to \RR$ by zeros. We note that the extended functions remain normalized.
Their respective Hamiltonian flows $\widetilde{\Phi}_t, \widetilde{\Psi}_t$ coincide with $\Phi_t, \Psi_t$ in $\TTs \setminus D_A$
and are constant in the complement.

Denote by $\gamma_\alpha, \gamma_\beta : S^1 \to (\TTs \setminus D_A) \subset \widetilde{\Sigma}$ the loops
$\gamma_\alpha (t) = (t, 0.5), \gamma_\beta (t) = (0.5, t)$.
They define homotopy classes $\alpha = [\gamma_\alpha], \beta = [\gamma_\beta] \in \pi_1 (\widetilde{\Sigma})$
whose intersection product $[\alpha] \cdot [\beta] = [pt]$ hence $\alpha$ and $\beta$ have topologically non-trivial intersection.

The composition of flows $\tilde{g}_A := \widetilde{\Psi}_{2A} \circ \widetilde{\Phi}_A \in \Ham(\widetilde{\Sigma})$ has exactly four
$1$-periodic points in the class $\alpha$. Indeed, we look for points
$(x, y) \in \TTs \setminus D_A$ where the horizontal velocity of $\Phi$ is $h'(y) = \frac{1}{A}$ and the vertical velocity of
$\Psi$ is $h'(x) = 0$. The four solutions are
\[
 \begin{array}{ll}
	  p_1 = \left(0, -\frac{1}{100 A}\right), &p_2 = \left(0.5, -\frac{1}{100 A}\right), \\
	  p_3 = \left(0, 0.5+\frac{1}{100 A}\right), &p_4 = \left(0.5, 0.5+\frac{1}{100 A}\right).
 \end{array}
\]
As $h'' (x), h'' (y) \neq 0$ at $p_i$, all four points are non-degenerate.

We choose $\gamma_\alpha$ to be a reference loop in the class $\alpha$ and pick a generic compatible almost complex structure
$J_t$. By direct computation, the actions of trajectories of $p_i$ under the the composition of flows
$\widetilde{\Psi}_{2A} \# \widetilde{\Phi}_A$ are
\[
  \begin{split}
	  \Ac (\gamma_{p_1}, w) &=  3A + \Delta_1 ,\\
	  \Ac (\gamma_{p_2}, w) &=  -A + \Delta_2 ,\\
	  \Ac (\gamma_{p_3}, w) &=   A + \Delta_3 ,\\
	  \Ac (\gamma_{p_4}, w) &= -3A + \Delta_4 .
  \end{split}
\]
Actions do not depend on cappings $w$, the error terms $|\Delta_i| < 1$.

Let $u$ be a non-trivial Floer cylinder connecting $(\gamma_{p_i}, w_i)$ to $(\gamma_{p_j}, w_j)$.
As $E(u) > 0$, the energy is bounded below by
\[
  E (u) = \Ac (\gamma_{p_i}, w_i) - \Ac (\gamma_{p_j}, w_j) = (\Ac_i + \Delta_i) - (\Ac_j + \Delta_j) > 2A - 2.
\]
By Corollary~\ref{Cor:stability} fixed points in the class $\alpha$
persist under $(2A-2)$-Hofer-small perturbations of $\tilde{g}_A$ in $\Ham(\widetilde{\Sigma})$.

A very similar argument applies to the class $\beta$.
$\tilde{g}_A$ has exactly four fixed points in this class, namely
\[
 \begin{array}{ll}
	  q_1 = \left(\frac{1}{200 A}, 0\right), &q_2 = \left(\frac{1}{200 A}, 0.5\right), \\
	  q_3 = \left(0.5-\frac{1}{200 A}, 0\right), &q_4 = \left(0.5-\frac{1}{200 A}, 0.5\right).
 \end{array}
\]
All four points are non-degenerate.

We pick $\gamma_\beta$ for a reference loop.
The actions of $q_i$ under the the composition of flows $\widetilde{\Psi}_{2A} \# \widetilde{\Phi}_A$ satisfy
\[
  \begin{split}
	  \Ac (\gamma_{q_1}, w) &=  3A + \Delta_1 ,\\
	  \Ac (\gamma_{q_2}, w) &=   A + \Delta_2 ,\\
	  \Ac (\gamma_{q_3}, w) &=  -A + \Delta_3 ,\\
	  \Ac (\gamma_{q_4}, w) &= -3A + \Delta_4 .
  \end{split}
\]
with the error term $|\Delta_i| < 1$.

As before, any non-trivial Floer cylinder connecting $(\gamma_{q_i}, w_i)$ to $(\gamma_{q_j}, w_j)$
has energy $E (u) > 2A - 2$. Hence by Corollary~\ref{Cor:stability} fixed points in the class $\beta$
persist under $(2A-2)$-Hofer-small perturbations of $\tilde{g}_A$.

That is, the obstruction of Proposition~\ref{P:obstruction} holds for any $h \in \Ham(\widetilde{\Sigma})$ in the
$(2A-2)$-Hofer neighborhood of $\tilde{g}_A$.
\begin{cor}
  $\dist_{H} \left(\tilde{g}_A, \Ham_{Aut} (\widetilde{\Sigma})\right) \geq 2A-2$.
\end{cor}

\begin{rmk}
  $\dist_{H} \left(\tilde{g}_A, \Ham_{Aut} (\widetilde{\Sigma})\right) \leq
  \dist_{H} \left(\tilde{g}_A, \widetilde{\Psi}_{2A}\right) \leq 2A$ hence the lower bound $2A-2$ is almost sharp,
  which is a rare case in Hofer's geometry. The additive constant $-2$ may be farther improved by fine tuning of the construction of $g_A$
  and computing precise action values for trajectories.
\end{rmk}


\subsection{A torus}
The argument is very similar to the previous subsection but we need to deal with an additional difficulty as
Floer cylinders may wrap around the torus. This may shorten the action gap between endpoint trajectories, the result is a smaller lower bound
for energy. We adjust the construction of $g_A \in \Ham (\TTs)$ in order to overcome this problem.

Fix a parameter $A > 2$. Recall that $F, P : \TTs \to \RR$ and the flows $\Phi_t, \Psi_t$ are supported in the complement of $D_A \subset \TTs$.
We modify the symplectic form in the disk $D_A$ so that the total area of the torus becomes $2$. Intuitively, this is equivalent
to puncturing $D_A$ at $q_0$ and gluing in a disk of area $1$.
We denote the resulting symplectic surface by $\TT$. That is, $\TT = \TTs$ as smooth manifolds, but the symplectic forms differ in $D_A$.
Consider $F, P$ as Hamiltonian functions $\TT \to \RR$, since the symplectic form was not altered in the supports
$\supp(F), \supp(P) \subset \TT \setminus D_A$, their induced flows are the same $\Phi_t, \Psi_t$.
We consider the same eggbeater map $g_A := \Psi_{2A} \circ \Phi_A$ as before.

Denote by $\gamma_\alpha, \gamma_\beta : S^1 \to \TT$ the loops
$\gamma_\alpha (t) = (t, 0.5), \gamma_\beta (t) = (0.5, t)$.
The homotopy classes $\alpha = [\gamma_\alpha], \beta = [\gamma_\beta] \in \pi_1 (\TT)$
have intersection product $[\alpha] \cdot [\beta] = [pt]$, hence $\alpha$ and $\beta$ have topologically non-trivial intersection.
By the same calculation as in the previous subsection, $g_A$ has exactly four $1$-periodic points in the class $\alpha$, namely,
$p_1, p_2, p_3, p_4$ and they are non-degenerate.

We pick $\gamma_\alpha$ as a reference loop in the class $\alpha$ and select a generic compatible almost complex structure
$J_t$. The actions of trajectories of $p_i$ under the the composition of flows $\Psi_{2A} \# \Phi_A$ are
\[
  \begin{split}
	  \Ac (\gamma_{p_1}, w) &=  3A + \Delta_1 + 2k_w ,\\
	  \Ac (\gamma_{p_2}, w) &=  -A + \Delta_2 + 2k_w ,\\
	  \Ac (\gamma_{p_3}, w) &=   A + \Delta_3 + 2k_w ,\\
	  \Ac (\gamma_{p_4}, w) &= -3A + \Delta_4 + 2k_w .
  \end{split}
\]
with the error term $|\Delta_i| < 1$ and $k_w \in \ZZ$ depends on the homotopy class of the capping $w$.
Very roughly, $k_w$ tracks how many times $w$ `wraps around' $\TT$, and it can be computed by the local degree
$k_w = \deg (w, q_0)$ of $w$ near $q_0$.

Let $u$ be a non-trivial Floer cylinder connecting $(\gamma_{p_i}, w_i)$ to $(\gamma_{p_j}, w_j)$.
If $q_0$ is not in the image of $u$, then both $w_i, w_j$ have the same wrapping number $k_i = k_j$.
Note that $p_i \neq p_j$ (otherwise the energy would be zero), therefore
\[
  E (u) = \Ac (\gamma_{p_i}, w_i) - \Ac (\gamma_{p_j}, w_j) = (\Ac_i + \Delta_i + 2k_i) - (\Ac_j + \Delta_j + 2k_j) > 2A - 2.
\]
If $q_0$ is in the image of $u$, the argument becomes a bit more subtle. First, $p_i \neq p_j$, as for a generic choice of data,
Floer cylinders cannot connect trajectories of the same Conley-Zehnder index.
(Different cappings do not affect the index since $T \TT$ is trivial.)
Suppose that
\[
  E (u) = \Ac (\gamma_{p_i}, w_i) - \Ac (\gamma_{p_j}, w_j) = (\Ac_i + \Delta_i + 2k_i) - (\Ac_j + \Delta_j + 2k_j) < A.
\]
Since the set $\{-3A, -A, A, 3A\}$ is $2A$-separated,
the `wrapping number' of $u$ given by $k_u = k_i - k_j = \deg (u, q_0)$ satisfies $|k_u| \geq \frac{A}{2}-1$.
We recall that both functions $F, P$ vanish in $D_A$, hence the Floer equation there coincides
with the equation of a $J$-holomorphic curve. The map $u$ is orientation-preserving and has constant degree (namely, $k_u$) in $D_A$.
Therefore $k_u > 0$ and
\[
  \begin{split}
	E (u) &= \int_{\RR \times S^1} \omega \left(\partial_s u, J_t (\partial_s u)\right) \dd t \dd s \geq
		  \int_{u^{-1} (D_A)} \omega \left(\partial_s u, J_t (\partial_s u)\right) \dd t \dd s = \\
		  &= \int_{u^{-1} (D_A)} \omega \left(\partial_s u, \partial_t u)\right) \dd t \dd s = k_u \cdot \Area (D_A) > k_u \geq \frac{A}{2}-1 .
  \end{split}
\]
That is, we obtain a lower bound $E (u) \geq \frac{A}{2}-1$.
Hence by Corollary~\ref{Cor:stability} fixed points in the class $\alpha$
persist under $\left(\frac{A}{2}-1\right)$-Hofer-small perturbations of $g_A$ in $\Ham(\TT)$.
The class $\beta$ is analyzed in a similar way, which results in the same lower bound $\left(\frac{A}{2}-1\right)$
for stability of its trajectories. Hence by of Proposition~\ref{P:obstruction}:
\begin{cor}
  $\dist_{H} \left(g_A, \Ham_{Aut} (\TT)\right) \geq \frac{A}{2}-1$.
\end{cor}
This finishes the proof of Theorem~\ref{T:main}.

\begin{rmk}
  The coefficient $\frac{1}{2}$ is not sharp.
  It can be improved to $(1 - \varepsilon)$ by increasing the relative area of the disk $D_A$.
\end{rmk}

%
%
%
%
%

\section{Discussion}\label{S:discuss}

Obstruction for a $g \in \Ham(\Sigma)$ to be autonomous which is used in this text
(namely, existence of two intersecting trajectories)
can be interpreted as an intermediate step between the one described in \cite{Po-Sh:egg-beater}
(existence of trajectories in a self-intersecting primitive homotopy class)
and that of \cite{Al-Me:braid} where the authors show that appearance of certain braids of trajectories implies
positive topological entropy.

That suggests the following question:
\begin{qs}
  Can one use similar methods to \cite{Al-Me:braid} to show that appearance of fixed points
  in both classes $\alpha$ and $\beta$ implies $h_\text{top} (g) > 0$?
\end{qs}
Preliminary considerations show that it is likely to be true, and in the case of a positive answer, the argument would imply that
$$\dist_H (g_A, \Ham_{\text{Ent}=0} (\Sigma)) \xrightarrow{A \to \infty} \infty .$$
A similar (and even stronger) statement is already established in
\cite{Ch-Me:Hof-entr} for surfaces of genus $g \geq 2$, our approach has a potential to extend it to genus $g \geq 1$.

If the answer is negative, one may consider a weaker question:
\begin{qs}
  Does the braid $\mathcal{B}$ formed by the eight fixed points of $\Fix_\alpha (g_A) \cup \Fix_\beta(g_A)$
  imply $h_\text{top} (g) > 0$?
\end{qs}
If the case of a positive answer, it is likely that the braid stability results of \cite{Al-Me:braid}
can be adapted to show that $\mathcal{B}$ persists under Hofer-perturbations of energy comparable to $A$.
Once again, this would imply $\dist_H (g_A, \Ham_{\text{Ent}=0} (\Sigma)) \to \infty$
as $A \to \infty$.

\medskip

We would like to add that while the eggbeater map $g_A$ was constructed in a rigid way with attention to small details,
the construction can be greatly flexified at the expense of a more delicate computation of energy bounds.
For example, one may avoid the perturbation of functions $F, P$ around $q_0$.
Or pick two independent parameters
$A, B > 1$ and let $g_{A, B} = \Psi_B \circ \Phi_A$.
One shows that $\dist_H (g_{A, B}, \Ham_{Aut}) \to \infty$ as $A, B \to \infty$, but in this setting the argument
relies on the Conley-Zehnder index in order to filter out certain potential cylinders of low energy.

\bibliography{bibliography}

\begin{thebibliography}{Cho}

\bibitem[AM]{Al-Me:braid}
M.~Alves and M.~Meiwes.
\newblock Braid stability and the {H}ofer metric.
\newblock Preprint (2021), can be found at arXiv:2112.11351v1.

\bibitem[BK]{Bra-Kha:C0-gap}
Michael Brandenbursky and Michael Khanevsky.
\newblock {$C^0$}-gap between entropy-zero {H}amiltonians and autonomous
  diffeomorphisms of surfaces.
\newblock Preprint, can be found at arXiv:2105.15038. To appear in the Israel
  J. Math.

\bibitem[Cho]{Ch:eggbeater}
Arnon Chor.
\newblock Eggbeater dynamics on symplectic surfaces of genus 2 and 3.
\newblock Preprint(2020). Can be found at https://arxiv.org/abs/2012.08930.

\bibitem[CM]{Ch-Me:Hof-entr}
Arnon Chor and Matthias Meiwes.
\newblock Hofer's geometry and topological entropy.
\newblock Preprint(2021). Can be found at https://arxiv.org/abs/2112.04955.

\bibitem[Hof]{Hof:TopProp}
H.~Hofer.
\newblock On the topological properties of symplectic maps.
\newblock {\em Proc. Royal Soc. Edinburgh}, 115(1-2):25--38, 1990.

\bibitem[Pol]{Po:geo-symp}
Leonid Polterovich.
\newblock {\em The Geometry of the Group of Symplectic Diffeomorphism}.
\newblock Lectures in Mathematics ETH Z{\"u}rich. Birkh{\"a}user, Basel, 2001.

\bibitem[PS]{Po-Sh:egg-beater}
Leonid Polterovich and Egor Shelukhin.
\newblock Autonomous {H}amiltonian flows, {H}ofer's geometry and persistence
  modules.
\newblock {\em Selecta Mathematica}, 22(1):227--296, Jan 2016.

\bibitem[Sal]{Sa:HF-lectures}
D.~Salamon.
\newblock Lectures on floer homology.
\newblock In {\em Symplectic geometry and topology (Park City, UT, 1997)},
  volume~7 of {\em IAS/Park City Math. Ser.}, pages 143--229. Amer. Math. Soc.,
  Providence, RI, 1999.

\bibitem[Sch]{Sc:action-spectrum}
M.~Schwarz.
\newblock On the action spectrum for closed symplectically aspherical
  manifolds.
\newblock {\em Pacific J. Math.}, 193(2):419--461, 2000.

\bibitem[Ush]{Ush:b-depth-appl}
M.~Usher.
\newblock Boundary depth in {F}loer theory and its applications to
  {H}amiltonian dynamics and coisotropic submanifolds.
\newblock {\em Israel J. Math.}, 184:1--57, 2011.

\end{thebibliography}

\end{document}